\input amstex
 \documentstyle{amsppt}
 \NoBlackBoxes

 \define \ov{\overline}
 
 \define \E{\Bbb E}
 \define \Q{\Bbb Q}
 \define \M{\frak M}
 \define \Z{\Bbb Z}
 \define \Hyp{{\Bbb H}^2}
 \define \sgn{\operatorname{sgn}}

 \topmatter

 \title
         non-positively curved graph manifolds  \\
         are virtually fibered over the circle
 \endtitle

 \author              P. Svetlov
 \endauthor

 \email
 \endemail

 \rightheadtext{Graph manifolds}

 \abstract
 In this note we prove that
 any closed graph manifold admitting a metric
 of non-positive sectional curvature (NPC-metric) has a finite
 cover, which fibers over the circle.

 An explicit criterion to have a finite cover, which fibers
 over the circle, is presented for the graph manifolds
 of certain class.
 \endabstract
 \endtopmatter

 \document
                    \title 0. Introduction \endtitle

 \subheading{0.1 Results} Which compact 3-manifolds
 admit Riemannian metrics of non-po\-si\-tive sectional
 curvature (NPC-metric)?
 It is known that if a 3-manifold $M$ admits a NPC-metric
 then it is irreducible
 and has infinite fundamental group.
 By the Thurston Geometrization Conjecture
 such a manifold is either hyperbolic
 or Seifert or Haken.
  In his paper \cite{3}, Leeb answered the question for
  Haken manifolds except for closed graph manifolds.
 (Recall that a graph manifold $M$ is a Haken manifold such that there are
 only Seifert components in the JSJ-decomposition of $M$).
 The last case was studied by Buyalo and Kobel'ski\u{i}
 in \cite{1,2} (see also Lemma 1.1).

 Another Thurston's conjecture
 claims that any {\it hyperbolic} 3-manifold
  is virtually fibered over the circle   (i.e.
 has a finite cover which is a surface bundle
 over the circle). Here we prove some ``similar"
 result.

 \proclaim{Theorem A}~If a closed graph manifold
 $M$ admits a NPC-metric
 then $M$  is virtually fibered over the circle.
 \endproclaim

 The assertion of theorem A
 is trivial for the graph manifolds
 with non-empty boundary. Indeed,
 {\it every} such manifold  is
 virtually fibered over the circle  \cite{9}
 as well as it admits a NPC-metric \cite{3}.
 On the other hand there is a lot of closed graph manifolds
 which do not admit NPC-metrics or are not virtually fibered.

 Let $\M$ be the class of closed orientable
 graph manifolds that are pasted from
 Seifert pieces with orientable base-orbifolds
 of negative Euler characteristic.
 Such Seifert manifolds admit a geometry
 modelled on $\Hyp\times\E^1$.

For each graph manifold $M\in\M$
 one can define a symmetric matrix $H_M$ (see 0.4),
 which consists
 of numerical invariants of $M$
 ($H_M$ is a generalization of
 a matrix introduced in \cite{2}).

 \proclaim{Theorem B}
 A graph manifold $M\in\M$
 is virtually fibered over the circle
 iff either $H_M$ has a negative eigenvalue or
 $H_M$ is supersingular (a matrix $A$ is called
 {\tt supersingular}
 if it annihilates a tuple with no zero entry \cite{5}).
 \endproclaim
\subheading{Remark} In the papers \cite{5,10}, two criteria of
virtual fibration for graph manifolds are proved. But application
of these criteria is too complicated, so they are implicit ones. In
the papers \cite{4,10}, two explicit {\it necessary} conditions of
virtual fibration for graph manifolds are presented.

 \subheading{0.2 The graph $\Gamma_M$}~
 Each graph manifold admitting a NPC-metric
 is finitely covered by a manifold of the class $\M$ \cite{1}.
 Such a covering manifold admits an (lifted) NPC-metric as well.
 In what follows we assume that all manifolds are in $\M$.

 Let $M$ be a graph manifold, let $\{M_v\}_{v\in V}$
 be the set of its Seifert pieces (we will call them
 {\it blocks}), and let
 ${\Cal T}=\cup_{e\in E}T_e$ be the JSJ-surface
 in $M$ (the minimal collection of incompressible tori).
 So this graph manifold $M=\cup_{\Cal T}M_v$
 can be obtained by pasting together
 the manifolds $M_v$ along its boundary tori.
 The collections $V$ and $E$ form a graph $\Gamma_M(V,E)$, which is dual
 to the JSJ-decomposition of $M$. Namely,
 two vertices $v,v'\in V$ (may be equal) are joined by an edge
 $e\in E$ if the blocks $M_v,\,M_{v'}\subset M$
 are (locally) separated by $T_e$. In this case we write $v'=e(v)$.
 By $\Gamma_M(V,W)$ we denote the oriented graph,
 which corresponded to $\Gamma_M(V,E)$, by $-w\in W$ we denote
 the opposite edge for $w\in W$, finally, by $\partial v$ we denote
 the set of edges (oriented or not) initiating at $v\in V$.

 \subheading {0.3 Invariants of graph manifolds}
 Choosing some orientation on $M$ we get an orientation on each
 maximal block $M_v,\,v\in V $ and hence
 on each its boundary torus
 $T_e\subset\partial M_v=\cup_{e\in\partial v}T_e$.
 The last orientation defines the canonical ``area" isomorphism
 $is_w:\Lambda^2H_1(T_e;\Q)\to\Q$, where the oriented
edge $w\in W$ corresponds to the orientation of $e\in E$ {\it
from} $v\in V$ (further we will write $w=(v,\,e)$).
 Let $a\wedge_w b=is_w(a\wedge b)\in \Q$ be the image of
 $a\wedge b\in\Lambda^2H_1(T_e;\Q)$ under this isomorphism.
 Note that
 $a\wedge_{-w} b=-a\wedge_{w} b$
 since the orientations
 on $T_e$ arising from $M_v$ and $M_{v'}$ are opposite.

 From now on we fix some orientations of the
 Seifert fibers in the blocks of $M$.
 The Seifert fibration of $M_v$ induces linear foliations
 of each boundary torus $T_e,\;{e\in\partial v}$.
 Let $f_w\in H_1(T_e,\Z),\;w=(v,\,e)$
 be the element representing the foliation of $T_{e}$.
 Note that the element $f_v=\iota_*(f_w)$,
 where the homomorphism $\iota_*:H_1(\partial M_v)\to H_1(M_v)$
 is induced by inclusion is independent
 of $w\in\partial v$ and presents any fiber of the Seifert fibration
 of $M_v$.

 We define the integer $b_w$ for each $v\in V$ and $w\in \partial v$
 by
 $$
 b_w=f_w\wedge_w f_{-w}.
 $$
 It is not difficult to see that $b_w=b_{-w}$.
 The integer $b_e=|b_w|,\;w=(v,\,e)$ is independent
 of the chosen orientations.

 The following lemma
 is useful to introduce another invariant.
 Its proof can be easily deduced from \cite{6} (see, also \cite{7}).

\proclaim{Lemma } Assume that for each $w\in\partial v$ we fix an
element $z_w\in H_1(T_e),\;w=(v,\,e)$ such that
$f_w\wedge_wz_w=1$. In this case we have
$$
\iota_*(z)=E(z)\cdot f_v,\qquad
\text{where}\;z=\bigoplus_{w\in\partial v}z_w,
$$
and $E(z)$ is a rational number (the Euler number).
\endproclaim
\subheading{Remark}
In paper \cite{6} the Euler number has the opposite sign.

Now we define  {\it the charge} $k_v$ of a block
 $M_v$ by
 $$
 k_v=E\left(\bigoplus_{\,\,w\in\partial v}
 \frac{f_{-w}}{b_w}\right).
 $$
 The charge $k_v$ of a block $M_v\subset M$
 is independent of the orientations of the Seifert fibers
 in blocks of M. It depends only on
 orientation of $M$.
 See, also \cite{1}.

 \subheading{0.4 The matrix $H_M$ \cite{2}}
 Consider the square matrix  $H_M=(h_{vv'})_{v,v'\in V}$
 over the graph $\Gamma_M(V,E)$:
 $$
 h_{vv'}=\left\{\matrix
 s(v)k_v&
 -\sum\limits_{e(v)=v}\frac{2}{b_e}&
 \text{if}\;\; v=v',
 \phantom{GGggGGGG}\\
 &
 -\sum\limits_{e(v)=v'}\frac{1}{b_e}&
 \phantom{G}\text{if}\;\;v\ne v'\;\;
 \text{and}\;\; k_vk_{v'}>0,\\
 &
 0 &
 \text{othe8ise},\phantom{GGGggGGG}\\
 \endmatrix\right.
 $$
 where the function $s:V\to \{0,\pm 1\}$ is constructed as follows.
 If $k_v=0$ for all $v\in V$ then $s(v)=0$ for all $v\in V$.
 Let $k_v\ne 0$ at least for one vertex $v\in V$.
 Vertices $v,v'\in V$
 of the graph $\Gamma_M$ are called equivalent ($v\sim v'$)
 if there exists  a path $v_0=v$, $v_1$, $v_2$, $\dots$, $v_n=v'$
 in $\Gamma_M$ such that
 $k_{v_i}\cdot k_{v_{i+1}}>0$ for each $i=0,\dots,\,n-1$.
 An edge $e\in E$ is called equivalent to a vertex $v\in V$ ($e\sim v$) if
 $k_v\ne 0$ and the vertices, incident
 to the edge are equivalent to $v\in V$.
 The factor graph  of $\Gamma_M(V,E)$ by the relation $\sim$ is called
 {\it the graph of signed components} and it is denoted by $G(U,E_0)$.
  For convenience we
 introduce the following notation
 $$
 p:\Gamma_M(V,E)\to G(U,E_0)=\Gamma_M(V,E)/\sim
 $$
 for the factor map.
 It is clear, that if $p\,(e)\in E_0$ then $e$
 joins nonequivalent vertices in $\Gamma_M(V,E)$.

 Now we are ready to define $s$.
If the graph $G(U, E_0)$ is not bipartite then we put
 $s(v)=0$ for each vertex  $v\in V$. Assume that $G(U,E_0)$ is
 a bipartite graph.
 Choose an orientation of $M$
 and a partition $U=P\cup N$ so that
 there exists a vertex $v\in P$ with $k_v>0$. In this case we put
 $s(v)=1$ if $v\in P$, and $s(v)=-1$ if $v\in N$.

 Note, that the properties of $H_M$ to be supersingular
 and to have a negative eigenvalue are independent of
 of representation of $H_M$ as a ``square table" (i.e. on
 an order of the vertices of $\Gamma_M$).
 On the other hand by reordering of vertices one can
 reduce this matrix to block-diagonal form
$H_M=\oplus_{u\in U}H_u$, where $H_u=\{h_{vv'}\,|\,v,v'\in p^{-1}(u)\}$.

 If the graph $G(U,E_0)$ is bipartite then this matrix
 coincides with the matrix $H_M$ from \cite{2}.
 If a graph manifold $M$
 has no block with zero charge then $G(U;E_0)$ is bipartite
 and $s(v)=\sgn k_v$. In this case
 $H_M$ can be represented as
 the matrix $-(P_-\oplus N)$ from the paper \cite{N1}.

 Throughout this paper we work in the piecewise linear category.

 \subheading{Acknowledgements}
 The paper is a part of the author's PhD dissertation
 and as such owes much to S. V. Buyalo.

 \head 1. Proof of theorem A
 \endhead
 Let  $\M_0\subset\M$  be the class of the closed oriented
 graph manifolds that are pasted from blocks of kind $F\times S^1$,
where $F$ is a surface of negative Euler characteristic.
 It is known that each manifold of the class $\M$ is finitely
 covered by a manifold from $\M_0$.
 Here we give results of \cite{1,2} which are used in what follows.

 \proclaim{Lemma 1.1} Let $M$ be a
 graph manifold of the class $\M_0$ and $\Gamma_M(V,E)$
 be its graph. The following conditions are equivalent.

 \roster
 \item $M$ admits a NPC-metric;
 \item the equation
$$
 k_va_v=\sum_{\overset w\in\partial v\to{
{\sssize w=(v,\,e)}}}
\frac{\gamma_ea_{e(v)}}{b_e}
 \tag{CE}
 $$
 over $\Gamma_M$
 with unknowns $\{a_v,\gamma_e\,|\,v\in V,\,e\in E\}$
 has a solution $\{a,\gamma\}$ such that
 $a_v>0$, $|\gamma_e|<1$ for each $v\in V$, $e\in E$;
 \item the matrix $H_M$ is either zero or has a negative eigenvalue.
 \endroster
 \endproclaim

 The equation \thetag{CE} is called {\it the Compatibility
 Equation over} $\Gamma_M$.

 \subheading{Remark}
 If \thetag{CE} has a solution then it has a rational one.
 Indeed, the conditions $a_v>0$, $|\gamma_e|<1$
 are ``open" ones and \thetag{CE} is a system of $|V|$
 quasilinear equations over $\Q$ in $|V|+|E|$ unknowns.
 Moreover, we can (and will) 
 assume that $a_v$ is integer for any $v\in V$.

 \medskip{\it Proof of theorem A.}~
 By the observationn above
 we can assume that $M\in\M_0$. If $M$ has a NPC-metric then
 the equation \thetag{CE}
 has a solution by lemma 1.1. Using this solution
 we are looking for a certain {\it immersed} surface
 in $M$ (Proposition 1.2).
 Then we show that this surface
 lifts as an {\it embedded} one in some finite cover of
 $M$ (Proposition 1.4).

 \proclaim{ Proposition 1.2}
 Let $M$ be a
 graph manifold of the class $\M$ and $\Gamma_M(V,E)$
 be its graph.

 If \thetag{CE} over $\Gamma_M(V,E)$
 has a rational solution $\{a,\gamma\}$
 such that $a_v>0$, $|\gamma_e|\le 1$ for each $v\in V$, $e\in E$
 then there exists a horizontal immersion
 $g:S\to M$ of an oriented surface with negative Euler characteristic
 to $M$ (an immersion of a surface to a graph manifold $M$
 is called {\tt horizontal}
 if it is transverse to the fibers of the Seifert
 fibered pieces of $M$).
 \endproclaim
 \subheading{Remark}
It is known \cite{8}, that the  homomorphism
 $g_*:\pi_1(S)\to \pi_1(M)$ induced by horizontal immersion
 is injective.

 \subheading{Remark}
 Here we give somewhat more general prove than is necessary
 for theorem A. We take $|\gamma_e|\le 1$ instead of $|\gamma_e|<1$
 since it useful for theorem B.

 \demo{Proof}
Let $\{a,\gamma\}$ be a rational solution of \thetag{CE}.
 We put $\gamma'_e=\sgn(b_w)\cdot\gamma_e$
 for $w=(v,\,e)$,~and consider the following classes in $H_1(T_e;\Q)$:
 $$
 c^{+}_{w}=\frac{1+\gamma'_e}{2b_w}
 \left(a_vf_{-w}+a_{e(v)}f_w\right)\,,\quad
 c^{-}_{w}=\frac{1-\gamma'_e}{2b_w}
 \left(a_vf_{-w}-a_{e(v)}f_w\right).
 $$
 We may assume (multiplying all $a_v$ by an appropriate positive integer)
 that the homological classes
 $c_{w}^{\pm}$ lie in $H_1(T_e;\Z)$.
 Since $-w=(e(v),\,e)$ we have  $c_{-w}^{+}=c_{w}^{+}$, 
 $c_{-w}^{-}=-c_{w}^{-}$.

 We choose a pair of
 nonoriented curves {\tt c}$_e^{\pm}$  on each JSJ-torus $T_e$ in $M$
  that realize (with some orientations)
 the classes $c_{w}^{\pm}$ respectively
 (if one of the classes is zero then there is no
 corresponding curve).

  Putting $c_w=c_w^++c_w^-$ it is easy to see that
 $$
 f_w\wedge_wc_{w}=a_v>0\quad
\text{is independent of}~w\in\partial v~\text{and}
\tag{1}
 $$
 $$
 \sum_{w\in\partial v}\frac{f_{-w}}{b_w}
 \wedge_{-w}c_{w}=
 \sum_{\overset w\in\partial v\to{
{\sssize w=(v,\,e)}}}\frac{\gamma'_ea_{e(v)}}{b_w}=
 k_v\cdot a_v\quad\text{by \thetag{CE}.}
 \tag{2}
 $$

 By $d(${\tt c}$)^{n}$ we denote $d$ copies of a curve obtained
 by going $n$ times around the curve $c$.
 It is proved in \cite{6, Lemma 3.1},
 that under two last conditions there exist
 integers $d_v,\,n_v$ so that for each integers $d,n$ which are divisible
 by $d_v$ and $n_v$ respectively there exists
 a {\it connected} horizontal immersed surface $S_v$ in $M_v$
 which spans the curves
 $d_v(${\tt c}$^{\pm}_{e})^{n_v},\,e\in\partial v$.
 (The converse assertion is also true:
 homological boundary of a horizontal surface, which
 is proper immersed in $M_v$, satisfies
 \thetag{1}, \thetag{2}). Taking appropriate integers $D$ and $N$ we
 get an immersed surface by fitting together the parts $S_v$ spanning the
 $D(${\tt c}$_e^{\pm})^N,\,e\in\partial v$.
 Such a surface can be realized by a horizontal immersion $g:S\to M$.
 \qed
 \enddemo

 In paper \cite{8}, Rubinstein and Wang  proved  a simple
 iff condition for a given
 horizontal surface in a graph manifold to lift
as an embedded surface in some finite cover
 (i.e be virtually embedded).
 Namely, let $S$ be a closed surface with negative Euler characteristic
 and let $g:S\to M\in \M$ be a horizontal immersion.
 Deforming $g$ slightly we may assume that $\Cal C=g^{-1}({\Cal T})$
 (recall that ${\Cal T}$ is a JSJ-surface of $M$)
 is a collection of disjoint simple closed curves in $S$.
 It is not difficult to see that each curve of this collection $\Cal C$
 is two-sided in $S$.
 Define a 1-cochain (a homomorphism of abelian groups)
 $s_g:C_1(S;\Z)\to \Q_+^*$ which takes values
 in the multiplicative group $\Q_+^*$ of positive rationals as follows.
 If $\kappa:[0,1]\to S$ is a singular 1-simplex
 and $k=[\kappa]\in C_1(S;\Z)$ is the corresponding 1-chain
 then $s_g(k)=1$ if the intersection $\kappa([0,1])\cap{\Cal C}$
 is empty, and
 $$
 s_g(k)=\left|\frac{f_w\wedge_wg_*[\,c\,]}%
 {f_{-w}\wedge_{-w}g_*[\,c\,]}\right|
 $$

 if the intersection $\kappa([0,1])\cap{\Cal C}$ is transversal
 and consists of one point
 on a curve $c\in {\Cal C}$,
  $g(c)\subset T_e$, and $w=(v,\,e)$.
 Here $[\,c\,]\in H_1(S;\Z)$ is the class of $c$
 (with some orientation),  $g\circ\kappa(0)\in M_{v}$
 and $g\circ\kappa(1)\in M_{e(v)}$.
 Such a map $s_g$ is extended (not uniquely) to entire group $C_1(S;\Z)$
 and the boundaries lie in the kernel of any extension.
 So we get some well defined cocycle {\tt s}$_g\in H^1(S;\Q_+^*)$.

 \proclaim{Lemma 1.3 \cite{8, theorem 2.3}} Suppose
 $M$ is an oriented graph manifold and $g:S\to M$ is a horizontal
 immersed surface.
 Then there exists a finite covering $p:\widetilde{M}\to M$
 and an embedding $\tilde{g}:S\to \widetilde{M}$
 such that $g=p\circ\tilde{g}$ if and only if {\tt s}$_g\equiv 1$.
 \endproclaim

 \proclaim{ Proposition 1.4}
Let $g:S\to M$ be the immersion, which was constructed
 in the proof of Proposition 1.2.
 Then we have {\tt s}$_g\equiv 1$.
  \endproclaim

 \demo{Proof}
 We choose a point $p_v\in g^{-1}(M_v)$ for each $v\in V$.
 Consider a JSJ-torus $T_e\subset M$
 that locally separates blocks $M_v$ and $M_{v'}$.
 The intersection $g(S)\cap T_e$
 consists on $D$ copies of $(\text{\tt c}^+_e)^N$
 and $D$ copies of $(\text{\tt c}^-_e)^N$.
 Let $k_i^+(w)$ (resp. $k_i^-(w)$) be a
 1-chain on $S$ such that
 $$\partial k_i^+(w)=[p_{v'}]-[p_{v}]
\quad (\text{resp.}\;\;\partial k_i^-(w)=[p_{v'}]-[p_{v}])$$
 and its support intersects the set ${\Cal C}\subset S$
 just one time
 by the pre-image of the $i$-th curve of $D(\text{\tt c}^+_e)^N$
 (resp. $D(\text{\tt c}^-_e)^N$), $i=1,\cdots,D$.
 Now we can compute the value of  $s_g$
 on the chains $k_i^{\pm}(w)$ ($v'=w(v)$):
$$
 s_g\Bigl(k_i^+(w)\Bigr)=\frac{f_w\wedge_{w}Nc_w^+}%
 {f_{-w}\wedge_{-w}Nc_{-w}^+}=\frac{a_v}{a_{v'}},\quad
 s_g\Bigl(k_i^-(w)\Bigr)=\frac{f_w\wedge_{w}Nc_w^-}%
 {f_{-w}\wedge_{-w}Nc_{-w}^-}=\frac{a_v}{a_{v'}},\quad
 $$

It is obvious, that $s_g$ is a coboundary, so
 {\tt s}$_g\equiv 1$.
 \qed

 To complete the proof of theorem A
 it remains to note that if a graph manifold $\widetilde{M}$
 contains an embedded  horizontal surface then
 either $\widetilde{M}$ itself or some its
 2-fold cover is fibered over the circle \cite{5}.
 \qed
 \enddemo

\head
              2.  Proof of theorem B
\endhead

 \proclaim{Proposition 2.1}
 Let $M$ be a
 graph manifold of the class $\M$ and $\Gamma_M(V,E)$
 be its oriented graph. The following
 two conditions are equivalent.
 \roster
 \item $M$ has a finite cover, which
 fibers over the circle;
 \item the compatibility equation
 over $\Gamma_M$  has a rational solution $\{a,\gamma\}$ such that
 $a_v>0$, $|\gamma_e|\le 1$ for each $v\in V$, $e\in E$.
\endroster
 \endproclaim

\demo{Proof}
We have already proved the implication $2\Rightarrow 1$ (propositions 1.2
and 1.4).

Assume that $M\in \M$ is virtually fibered.
Let $g:S\to M$ be a virtual embedded horizontal surface 
(i.e. g(S) is the image of a fiber
 under the covering) and
$c_i,~{i\in I}$ be curves of the collection ${\Cal C}=g^{-1}({\Cal
T})$ on $S$ (${\Cal T}$, as usually, denotes the JSJ-surface in
$M$). The set $S\setminus{\Cal C}$ is a disjoint union of connected
components $S_{\alpha},~{\alpha\in A}$. Choosing orientation of $M$
and orientations of the Seifert fibers in its blocks provide
orientations of images $g(S_{\alpha})$ and hence ones of each component
$S_{{\alpha}}$.
 Elements of the collections $A$ and $I$ are
 vertices and edges respectively of a graph $\Gamma_S(A,I)$, which is dual
 to the decomposition of $S=\cup_{\Cal C}S_{\alpha}$ along the
 curves of $\Cal C$. Namely,
 two vertices $\alpha,\,\beta\in A$ (may be equal) are joined by an edge
 $i\in I$ if the components $S_{{\alpha}}$, $S_{{\beta}}$
 are (locally) separated by $c_i$. In this case we write $\beta=i(\alpha)$.
 By $\Gamma_S(A,L)$ we denote the oriented graph,
 which corresponded to $\Gamma_S(A,I)$. By $-l\in L$ we denote
 the opposite edge for $l\in L$. Finally, by $\partial \alpha$ we denote
 the set of edges (oriented or not) initiating at $\alpha\in V$.

 By $r:\Gamma_S(A,I)\to\Gamma_M(V,E)$
we will denote the natural map that preserves the incidence
relation and maps the vertices to vertices
 and the edges to edges so that $g(S_\alpha)\subset M_{r(\alpha)}$
and $g(c_i)\subset T_{r(i)}$. It is clear that the map $r$ is lifted
to a map between oriented graphs (we beep use the same notation for this
lifted map).

 Now we fix a continuous embedding $\tau:\Gamma_S(A,I)\to S$ such that
 $\tau(\alpha)\in S_{\alpha}$ and the intersection
 $\tau(i)\cap{\Cal C}$ is transversal and consists of one point
 on a curve $c_i$ for each $i\in I$. 
 In this situation one can consider the oriented edges of
 $\Gamma_S$ as 1-simplexes,  and their images
 as 1-simplexes on  $S$.
 If $S$ is virtually embedded then
 $\tau^*\text{\tt s}_g$ is the trivial element of $H^1(\Gamma_S(A,I);\Q^*_+)$
 by lemma 1.3, i.e. the 1-cochain $\tau^*s_g$
 is a coboundary. So there exists a function $j:A\to\Q^*_+$ such that
 for each oriented edge $l\in L$ we have
 $\tau^*s_g(l)=j_{\sssize\alpha}/j_{\sssize l(\alpha)}$.
On   the other hand
$$
\tau^*s_g(l)=s_g(\tau_*l)=\frac{f_{r(l)}\wedge_{r(l)}g_*[c_i]_{\alpha}}%
{f_{r(-l)}\wedge_{r(-l)}g_*[c_i]_{l(\alpha)}},
$$
where $[c_i]_{\alpha}\in H_1(S)$
is the homological class of
$c_i$ which is oriented as a part of
$\partial \ov{S}_{{\alpha}}$
and $l=(\alpha,\,i)$.
Comparing the obtained equalities for
$\tau^*s_g(l)$ we get
$$
f_{r(l)}\wedge_{r(l)}g_*[c_i]_{\sssize\alpha}=\phi_i\cdot j_{\sssize\alpha},
$$
where $\phi:I\to\Q^*_+$ is a function.

Now we apply lemma 3.1 from \cite{6} (see the end of the proof
of proposition 2.1) to the proper immersion
$g|_{\ov{S}_\alpha}:(\ov{S}_\alpha,\partial\ov{S}_\alpha)\to
(M_{r(\alpha)},\partial M_{r(\alpha)})$.
Consider the element $c_w^\alpha\in H_1(T_e;\Z)$, $w=(r(\alpha),\,e)$
that represents the  boundary of $g(\ov{S}_\alpha)$
on $T_e$:
$$
c_w^\alpha=\sum_{\sssize\beta\in w^r(\alpha)}
\sum_{i\in I_{\sssize\alpha e\beta}}g_*[c_i]_\alpha,
$$
where $w^r(\alpha)=\{\beta\in A\,|\;r(\beta)=w\left(r(\alpha)\right)\}$
and $I_{\alpha e\beta}=\{i\in I\,|\,r(i)=e,\,i(\alpha)=\beta\}$.
Then one can rewrite \thetag{1}, \thetag{2}
as
$$
a_\alpha=f_{w}\wedge_{w}c_w^\alpha=\left( \sum_{\sssize\beta\in
w^r(\alpha)} \sum_{i\in I_{\sssize\alpha
|w|\beta}}\phi_i\right)\cdot j_\alpha ~\text{is independent
of}~w\in\partial r(\alpha), \tag{$1'$}
$$
$$
k_{r(\alpha)}a_\alpha =\sum_{w\in\partial r(\alpha)}\frac{1}{b_w}
\sum_{\sssize\beta\in w^r(\alpha)}
\sum_{i\in I_{\sssize\alpha |w|\beta}}f_{-w}\wedge_{-w}g_*[c_i]_\alpha,
\tag{$2'$}
$$
where $|w|$ denotes the nonoriented edge of $\Gamma_M$
corresponding to $w\in W$.
If a curve $c_i$ separates components $S_\alpha$ and $S_\beta$
then $[c_i]_\alpha=\epsilon_i[c_i]_\beta$, $\epsilon_i=\pm 1$,
hence
$$
\sum_{i\in I_{\sssize \alpha |w|\beta}}f_{-w}\wedge_{-w}g_*[c_i]_\alpha=
\sum_{i\in I_{\sssize \alpha |w|\beta}}
\epsilon_if_{-w}\wedge_{-w}g_*[c_i]_\beta=\left(
\sum_{i\in I_{\sssize \alpha |w|\beta}}\epsilon_i\phi_i\right)\cdot j_\beta.
$$
Put
$$
\Phi^e_{\sssize\alpha\beta}=\sum_{i\in I_{\alpha e\beta}}\phi_i,\quad
\gamma^e_{\sssize\alpha\beta}=
\frac{1}{\Phi^e_{\sssize\alpha\beta}}
\sum_{i\in I_{\sssize \alpha e\beta}}\epsilon_i\phi_i\;.
$$
We assume that $\Phi^e_{\sssize\alpha\beta}=\gamma^e_{\sssize\alpha\beta}=0$
if $I_{\sssize \alpha e\beta}=\emptyset$.
Now we can rewrite the equalities \thetag{$1'$} and  \thetag{$2'$}
as
$$
a_\alpha=\left(\sum_{\beta\in w^r(\alpha)}
\Phi_{\sssize\alpha\beta}^{|w|}\right)\cdot j_\alpha\,
\qquad\text{is independent of}~w\in\partial r(\alpha),
\tag{$1''$}
$$
$$
k_{r(\alpha)}a_\alpha =\sum_{w\in\partial r(\alpha)}\frac{1}{b_w}
\sum_{\sssize\beta\in w^r(\alpha)}
\Phi_{\sssize\alpha \beta}^{|w|}j_\beta.
\tag{$2''$}
$$

Let $\Phi_w$ be the rectangular matrix with entries
$\left(\gamma^{|w|}_{\sssize\alpha\beta}
\cdot\Phi^{|w|}_{\sssize\alpha\beta}\right)$,
where the subscripts $\alpha,\;\beta$ range over 
$r^{-1}(v),\;r^{-1}\left(w(v)\right)$ respectively (here
$v$ is the initiating vertex of $w\in W$).
Note, that $\Phi_{-w}=\Phi_w^t$, where $t$ denotes
the transposition. Let
$\Phi_v$ be the diagonal matrix with entries
$\left(a_{\sssize\alpha}/j_{\sssize\alpha}\right)$ and let
$J_v$ be the tuple with entries
$(j_\alpha)$, where the subscript
$\alpha$ ranges over $r^{-1}(v)$.

Now we draw the equality \thetag{$2''$}
in new notations as a vector one:
$$
k_v \Phi_vJ_v=\sum_{w\in \partial v} \frac{1}{b_w}\Phi _wJ_{w(v)}.
\tag{$2'''$}
$$

Put
$$
x_v=\sqrt{J_v^t\Phi_vJ_v}>0,
\qquad\gamma_w=\sgn(b_w)\frac{J_v^t\Phi_wJ_{w(v)}}
{x_{w(v)}x_v},
$$
then by \thetag{$2'''$} we have
  $$
  k_vx_v=\sum_{w\in\partial v}\frac{\gamma_wx_{w(v)}}{|b_{w}|}.
  $$

\proclaim{Lemma 2.2} $|\gamma_w|\le 1$.
\endproclaim
{\it Proof.}~
If $A=(a_{ik})$ is a $m\times n$-matrix,
$x=(x_i)$ is a $m$-tuple,
$y=(y_k)$ is a $n$-tuple, and all numbers
$a_{ik},\;x_i,\;y_k$ are non-negative then
$$
\multline
x^tAy=\sum_{i,k}x_ia_{ik}y_k=
\sum_ix_i\left(\sum_k\sqrt{a_{ik}}\cdot\sqrt{a_{ik}}y_k\right)\le\\
\sum_i\left(x_i\sqrt{\sum_ka_{ik}}\right)\cdot
\sqrt{\sum_ka_{ik}y_k^2}\le
\sqrt{\sum_i\left(\sum_ka_{ik}\right)x_i^2}\cdot
\sqrt{\sum_k\left(\sum_ia_{ik}\right)y^2_k}
\endmultline
$$
by Cauchy-Bunyakovskii inequality applyied twice.

Since
$|\gamma^e_{\sssize\alpha\beta}|\le 1$ for all subscripts, we
have $\left|J^t_v\Phi_wJ_{w(v)}\right|\le J^t_vA_wJ_{w(v)}$, where
$A_w$ is rectangular matrix with entries
$\left(\Phi^{|w|}_{\sssize\alpha\beta}\right)$.
Previous sequence of inequalities, applyed to the matrix $A_w$
and to the tuples $J_v,\;J_{w(v)}$, gives
$$
J^t_vA_wJ_{w(v)}\le
\sqrt{J^t_v\Phi_vJ_v\cdot J^t_{w(v)}\Phi_{w(v)}J_{w(v)}}.
$$
This completes the proof of the lemma. \qed

To complete the proof of proposition 2.1
it remains to note that $\gamma_{-w}=\gamma_w$ since
$\Phi_w^t=\Phi_{-w}$.
In such a way the collection $\{x,\gamma\}$ give the required
solution of \thetag{CE}.
\qed
\enddemo

\demo{The proof of theorem B}
If
$H_M$
has a negative eigenvalue then
$M$
admits a NPC-metric (by Lemma 1.1).
So
$M$
is covered by a surface bundle (by theorem A).
Therefore to prove the ``if" part of theorem B
it remains to assume that
$G(U,E_0)$
is bipartite,
$H_M\ge 0$
(i.e.
$x^tH_Mx\ge 0$
for each tuple
$x$),
and
$H_M$
is supersingular. By supersingularity, there exists a tuple
$l$
with no zero entry such that
$H_Ml=0$.
We claim that all entries of
$l$
have the same sign.
Indeed, for each tuple
$x$
we can write the following well known identity (see, e.g. \cite{6})
$$
x^tH_Mx=-\frac{1}{2}\sum_{v,v'\in V}h_{vv'}l_vl_{v'}
\left(\frac{x_{v}}{l_{v}}-\frac{x_{v'}}{l_{v'}}\right)^2
=\frac{1}{2}\sum_{v,v'\in V}|h_{vv'}|\,l_vl_{v'}
\left(\frac{x_{v}}{l_{v}}-\frac{x_{v'}}{l_{v'}}\right)^2.
$$
We put
$x_v=|l_v|$
for each
$v\in V$
and then write
$$
x^tH_Mx=\frac{1}{2}\sum_{v,v'\in V}|h_{vv'}|\,l_vl_{v'}
\left(\sgn{l_{v}}-\sgn{l_{v'}}\right)^2.
$$
The left hand side of the equality is greater than or equal to zero
since $H_M\ge 0$ and the right hand side of the equality
is less than or equal to zero.
Therefore we have
$\sgn{l_{v}}=\sgn{l_{v'}}$
for each
$v,v'\in V$.

Let $a_v=|l_v|$ for each $v\in V$ be a ``positive"
tuple and
$$\gamma_e=
\left\{
\matrix
\phantom{-}1& \text{if}&p\,(e)\in P,\\
-1&           \text{if}&p\,(e)\in N,\\
\phantom{-}0& \text{if}&\,p\,(e)\in E_0.
\endmatrix
\right.
$$
For notations $P,\,N,\,E_0$ see 0.4.
It is not difficult to see
that $\{a,\gamma\}$
is a solution of \thetag{CE} over
$\Gamma_M$. Therefore $M$ is virtually fibered by proposition 2.1.

Conversely, by proposition 2.1 we can start from a solution
$\{a,\gamma\}$ of CE over
$\Gamma_M$. If
$|\gamma_e|<1$ for each
$e\in E$ then there is nothing to prove since in this case
(and only in this case if
$H_M\not\equiv 0$)
$H_M$ has a negative eigenvalue by lemma 1.1.
So we assume that $H_M\ge 0$ and $H_M\not\equiv 0$.
Let
$H_M=\oplus_{u\in U}H_u$ be the block-diagonal expansion
corresponding to the graph $G(U,E_0)$ and let
$a^u$ be the sub-tuple of $a$ corresponding to
$u\in U$. It is clear, that $H_M\ge 0$ if and only if
$H_u\ge 0$ for any
$u\in U$.
Putting
$V_u=p^{-1}(u)\cap V$ and
$E_u=p^{-1}(u)\cap E$
we can write
$$
(a^u)^tH_u\,a^u=\sgn u\sum_{v\in V_u}k_va_v^2-
2\sum_{e\in E_u}\frac{a_{e^+}a_{e^-}}{b_e},
$$
where $\sgn u=1$ if
$u\in P$, $\sgn u=-1$ if
$u\in N$, and  $e^{\pm}$ are the vertices that are joined by
$e$. Using CE we can rewrite
$$
\multline
(a^u)^tH_u\,a^u=\sgn u\sum_{v\in V_u}
\sum_{e\in \partial v}\gamma_e\frac{a_{e^+}a_{e^-}}{b_e}
-2\sum_{e\in E_u}\frac{a_{e^+}a_{e^-}}{b_e}=\\
\sgn u\sum_{v\in V_u}
\left(\sum_{e\in \partial v\cap V_u}\gamma_e\frac{a_{e^+}a_{e^-}}{b_e}
+\sum_{e\in \partial v\cap E_0}\gamma_e\frac{a_{e^+}a_{e^-}}{b_e}\right)
-2\sum_{e\in E_u}\frac{a_{e^+}a_{e^-}}{b_e}=\\
\sgn u\sum_{e\in E_0^u}\gamma_e\frac{a_{e^+}a_{e^-}}{b_e}
-2\sum_{e\in E_u}(1-\sgn u\cdot\gamma_e)\frac{a_{e^+}a_{e^-}}{b_e},
\endmultline
$$
where $E_0^u=\{e\in E_0\,|\,
e~\text{has one of its end in}~V_u\}$.
Summing the last equalities for
$u\in U$ we get
$$
a^tH_M\,a=\sum_{u\in U}(a^u)^tH_u\,a^u=
-2\sum_{e\in E\setminus E_0}(1-\sgn e\cdot\gamma_e)
\frac{a_{e^+}a_{e^-}}{b_e},
$$
where $\sgn e=1$ if
$p\,(e)\in P$ and $\sgn e=-1$ othewise. But $a^tH_Ma\ge 0$,
therefore
$\gamma_e=\sgn e$ for each $e\in E\setminus E_0$ and
$a^tH_Ma=0$. Since this matrix $H_M$ is symmetric and
 positive semidefined, we have $H_Ma=0$.
 But $a_v>0$ for each
$v\in V$, so
$H_M$ is supersingular.
\qed
\enddemo

\medskip
e-mail: {\tt svetlov\@pdmi.ras.ru}

 \end{document}